\documentclass[psamsfonts,fceqn,leqno]{amsart}
  \usepackage{amsthm}
  \usepackage{amsmath}
  \usepackage{amssymb}
  \usepackage{color}
  \usepackage{mathrsfs}
  
  \newtheorem{theorem}{Theorem}[section]
  
  \newtheorem{proposition}[theorem]{Proposition}
  
  \newtheorem{question}[theorem]{Question}
  \newtheorem{problem}[theorem]{Problem}
  \theoremstyle{definition}

  \newtheorem{example}[theorem]{Example}

  \numberwithin{equation}{section}

  \title[Three open problems on the Wijsman topology]
  {Three open problems on the Wijsman\\ topology$^\ast$}
  \author[Jiling Cao]{Jiling Cao}
  \address{School of Engineering, Computer and Mathematical Sciences, Auckland University of Technology, 
  Private Bag 92006, Auckland 1142, New Zealand}
  \email{jiling.cao@aut.ac.nz}
  \thanks{\hspace{-1.66em} $^\ast$This article was based on the author's talk at the 16$^{\rm th}$
  Galway Topology Colloquium, Galway, Ireland, 8--10 July 2013. The author wishes to thank the School 
  of Mathematics, Statistics and Applied Mathematics at NUI Galway for its hospitality during his visit
  in July 2013.}
  \date{}

\begin{document}

  \begin{abstract}
  Since it first emerged in Wijsman's seminal work \cite{wijsman:66}, the Wijsman topology has been intensively 
  studied in the past 50 years. In particular, topological properties of Wijsman hyperspaces, relationships 
  between the Wijsman topology and other hyperspace topologies, and applications of the Wijsman topology in 
  analysis have been explored. However, there are still several fundamental open problems on this topology. 
  In this article, the author gives a brief survey on these problems and some up-to-date partial solutions.
  \end{abstract}

  \maketitle

  \section{Introduction} \label{sec:intro}
  
  Let $(X,d)$ be a metric space, and $2^X$ be the collection of all non-empty closed subsets of $X$. There
  are many ways to equip a topology $\mathscr T$ on $2^X$ such that $(X, {\mathscr T}(d))$ is embedded into 
  $\left(2^X, {\mathscr T}\right)$ as a closed subspace via the mapping $x \mapsto \{x \}$, where $\mathscr
  T(d)$ is the topology on $X$ induced by the metric $d$. A unified approach to topologize $2^X$ discussed 
  in the monograph \cite{beer:93} is using distance functionals $d(\cdot, A): X \to \mathbb R^+$ for sets 
  $A \in 2^X$, where $d(\cdot, A)$ is defined such that for each $x\in X$,
  \[
  d(x, A) =\inf \{d(x, a); a \in A\}.
  \] 
  The weakest topology on $2^X$ such that all functionals $d(\cdot, A)$, $A \in 2^X$, are continuous is 
  called the \emph{Wijsman topology} and denoted by $\mathscr T_{W(d)}$. Correspondingly, we call $\left(
  2^X, \mathscr T_{W(d)}\right)$ the \emph{Wijsamn hyperspace} of $(X,d)$. The Wijsman topology is formally 
  introduced in Lechicki and Levi \cite{lechicki-levi:87}, but it can be tracked back to the seminal work
  of Wijsman \cite{wijsman:66},  where R.~A. Wijsman considered a mode of convergence for sequences of 
  closed sets when he studied some optimum properties of sequential probability ratio test in 1960's.
  
  \medskip
  The Wijsman topology is closely related to the other two well studied hyerspace topologies on $2^X$: the 
  \emph{Hausdorff} \emph{metric topology} $\mathscr T_{H(d)}$ and the \emph{Vietoris topology} $\mathscr 
  T_V$. Firstly, it is easily verified that the (extended) Hausdorff distance $H_d(A,B)$ between 
  two members $A, B \in 2^X$ can be defined as
  \[
  H_d(A,B) = \sup\{|d(x, A) -d(x, B)|: x\in X \}.
  \]
  Thus, the Hausdorff metric topology is just the topology of uniform convergence on $2^X$ under the 
  identification $A \leftrightarrow d(\cdot, A)$, while the Wijsman topology is the topology of pointwise 
  convergence on $2^X$ under the same identification. Secondly, it is known that the Vietoris topology on 
  $2^X$ is given by
  \[
  \mathscr T_V=\sup \left\{\mathscr T_{W(\varrho)}: \varrho \mbox{ is a metric on $X$ equivalent to $d$} 
  \right\}.
  \]
  The reader can find more details on relationships among these and other hyperspace topologies in 
  \cite{beer:93}. However, it turns out that the Wijsman topology is more intractable than both of the 
  Hasudorff metric and Vietoris topologies.
  
  \medskip
  In the past 50 years, the Wijsman topology and convergence under this property has been intensively 
  studied. For example, possible extensions of Wijsman's results in infinite-dimensional Banach spaces 
  have been considered, refer to \cite{beer:94}. In addition, various properties of the Wijsman topology, 
  conditions for the coincidence of the Wijsman topology with other hyperspace topologies and function 
  spaces equipped with the Wijsman topology have also been investigated, refer to  \cite{beer:93}, 
  \cite{beer:94} and \cite{diMM:98}. The Wijsman topology in the context of quasi-metric spaces has been 
  studied in \cite{cao-lopez:12} and \cite{lopez-romaguera:03}. In the study of properties of the Wijsman 
  topology, several techniques including game-theoretic approaches, special embeddings and construction 
  of special power spaces etc., have been employed. Despite of these efforts, several fundamental 
  problems still remain unsolved. 
  
  \medskip
  In this short article, the author will discuss three basic open problems on the Wjsman topology and some 
  open question associated with these problems, studied by him and his co-authors in the past 10 years. 
  The first problem is when the Wijsman topology induced by a metric is normal, the second problem is 
  whether the Wijsman topology induced by a completely metrizable metric has some completeness-type 
  properties, and the third one is to find characterization as when the Wijsman topology has the Baire 
  property in terms of some properties of the underlying metric space. For each of these problems, the 
  author will give a brief survey on its background and then also some up-to-date partial solutions. The 
  reader can find more details of these problems and some of their associated questions from the listed 
  references, particularly from \cite{cao:10}, \cite{cao-junnila:10}, \cite{cao-junnila:2014}, 
  \cite{cao-junnila:12} and \cite{cao-tomita:07}, respectively. Note that the author has no intention 
  to give a comprehensive and up-to-date survey on the Wijsman topology. For undefined notation and 
  background knowledge on the Wijsman topology, refer to \cite{beer:93}, \cite{beer:94} and \cite{diMM:98}.

  \section{The Normality Problem and Associated Questions}~\label{sec:normality}

  First of all, note that all Wijsman topologies are weak topologies, and thus they are Tychonoff. 
  However, not all Wijsman topologies are normal. To see this, let $d$ be the $0$-$1$ metric on a 
  nonempty set $X$ with $|X|=\aleph_1$. As observed in Remark 3.1 of \cite{chaber-pol:02}, $\left(2^X, 
  \mathscr T_{W(d)}\right)$ is homeomorphic to $\{0,1\}^{\aleph_1}\setminus\{ {\bf 0}\}$, where $\{ 0, 
  1\}$ is equipped with the discrete topology and ${\bf 0}$ is the constant function with value $0$. It 
  is known that $\{0, 1\}^{\aleph_1}\setminus \{ {\bf 0}\}$ is not normal, and consequently, $\left(2^X, 
  \mathscr T_{W(d)}\right)$ is not a normal space. Thus, the following natural problem arises.

  \begin{problem}[\cite{cao-junnila:2014}] \label{prob:normal}
  Let $(X,d)$ be a metric space. When is the Wijsman hyperspace $\left(2^X, \mathscr T_{W(d)}\right)$ 
  a normal space?
  \end{problem}

  A classical result of Lechicki and Levi in \cite{lechicki-levi:87} states that for a metric space 
  $(X, d)$, $\left(2^X, \mathscr T_{W(d)}\right)$ is metrizable if and only if $(X,d)$ is separable.
  As an immediate consequence, if $(X,d)$ is separable, then $\left(2^X, \mathscr T_{W(d)}\right)$ is 
  normal. However, it is unclear if the converse holds. Indeed, the following problem was posed by 
  Di Maio and Meccariello in 1998.

  \begin{problem}[\cite{diMM:98}]\label{prob:dimm}
  It is known that if $(X,d)$ is a separable metric space, then $\left(2^X,\mathscr T_{W(d)}\right)$ 
  is metrizable and so paracompact and normal. Is the opposite true? Is $\left(2^X, \mathscr T_{W(d)}
  \right)$ normal if, and only if, $\left(2^X, \mathscr T_{W(d)}\right)$ is metrizable?
  \end{problem}
  
  Regarding Problems \ref{prob:normal} and \ref{prob:dimm}, it was conjectured that for a metric 
  space $(X,d)$, if $\left(2^X, \mathscr T_{W(d)}\right)$ is normal, then $(X,d)$ is separable. In other
  words, $\left(2^X, \mathscr T_{W(d)}\right)$ is non-normal for a non-separable metric $(X,d)$. If this
  conjecture is true, then the answer to Problem \ref{prob:dimm} is affirmative, and also a
  characterization of normality of $\left(2^X, \mathscr T_{W(d)}\right)$ is derived and thus a solution
  to Problem \ref{prob:normal} is obtained. Below, I shall summarize some progress in this direction.
  
  \medskip  
  For a metric space $(X,d)$, define ${\rm nlc}(X)$ by
  \[
  {\rm nlc}(X) = \{x \in X:\ x \mbox{ has no compact neighbourhood
  in } X\}.
  \]
  The following embedding theorem was proved by Chaber and Pol in 2002.

  \begin{theorem}[\cite{chaber-pol:02}] \label{thm:chaber-pol02}
  Let $X$ be a metrizable space such that $w({\rm ncl}(X))= 2^{\aleph_0}$. Then for any compatible 
  metric $d$, ${\mathbb N}^{2^{\aleph_0}}$ embeds as a closed subspace in $\left(2^X, \mathscr T_{W(d)}
  \right)$. In particular, $\left(2^X, \mathscr T_{W(d)}\right)$ contains a closed copy of $\mathbb Q$.
  \end{theorem}

  Following the proof of Theorem \ref{thm:chaber-pol02}, if ${\rm nlc}(X)$ is a non-separable subspace 
  of $(X,d)$, $\left(2^X, \mathscr T_{W(d)}\right)$ contains a closed copy of ${\mathbb  N}^{\aleph_1}$. 
  This implies that $\left(2^X, \mathscr T_{W(d)}\right)$ is non-normal if ${\rm nlc}(X)$ is a 
  non-separable subspace of $(X,d)$, since ${\mathbb N}^{\aleph_1}$ is non-normal. Particularly, if 
  $(X,\|\cdot\|)$ is a non-separable normed linear space and $d$ is the metric on $X$ induced by 
  $\|\cdot \|$, then $\left(2^X, \mathscr T_{W(d)}\right)$ is non-normal. As a consequence, we derive a 
  partial answer to Problems \ref{prob:normal} and \ref{prob:dimm}, due to Hol\'{a} and Novotn\'{y} in 
  2013.

  \begin{theorem}[\cite{hola-novotny:13}]
  Let $(X, \|\cdot\|)$ be a normed linear space, and let $d$ be the metric on $X$ induced by 
  $\|\cdot \|$. Then $\left(2^X, \mathscr T_{W(d)}\right)$ is normal if and only if it is metrizable.
  \end{theorem}

  Suppose $(X,d)$ is a non-separable metric space. Then $(X,d)$ contains an $\varepsilon$-discrete 
  subspace $Y$ of size $\aleph_1$ for some $\varepsilon >0$. It follows that $2^Y$ is closed subspace 
  of $\left(2^X, \mathscr T_{W(d)}\right)$. Thus, if one can show that $2^Y$ is a non-normal subspace 
  of $\left(2^X, \mathscr T_{W(d)}\right)$, then Problems \ref{prob:normal} and \ref{prob:dimm} would 
  be solved. Indeed, to the author's knowledge, the following question is still unsolved.

  \begin{question} \label{ques:uniformdiscrete}
  Let $(X,d)$ be a uniformly discrete and non-separable metric space. Must $\left(2^X, \mathscr T_{W(d)}
  \right)$ be non-normal?
  \end{question}
  
  Inspirited by the work of Keesling in \cite{keesling:70a}, Cao and Junnila \cite{cao-junnila:2014}
  explored the possibility whether the non-normal space $\omega_1 \times (\omega_1 +1)$ is embeddable
  into the Wijsman hyperspace of a non-separable metric space $(X,d)$, and they obtained the following 
  result.

  \begin{proposition}[\cite{cao-junnila:2014}] \label{prop:subspace}
  Let $(X,d)$ be a non-separable metric space. Then the subspace $2^X\setminus \{X\}$ of $\left(2^X, 
  \mathscr T_{W(d)}\right)$ contains a closed copy of the space $\omega_1\times (\omega_1+1)$.
  \end{proposition}
  
  Applying Proposition \ref{prop:subspace} and the classical result of Lechicki and Levi in 
  \cite{lechicki-levi:87} on metrizability of Wijsman hyperspaces, Cao and Junnila were able to
  derive the following partial answer to Problems \ref{prob:normal} and \ref{prob:dimm}.  
  
  \begin{theorem}[\cite{cao-junnila:2014}] \label{thm:hnormal}
  Let $(X,d)$ be a metric space. The following are equivalent.
  \begin{itemize}
  \item[(1)] $\left(2^X, \mathscr T_{W(d)}\right)$ is metrizable.
  \item[(2)] $\left(2^X, \mathscr T_{W(d)}\right)$ is hereditarily normal.
  \item[(3)] $2^X\setminus \{X\}$ is a normal subspace of $\left(2^X, \mathscr T_{W(d)}\right)$.
  \end{itemize}
  \end{theorem}
  
  Although the techniques of embeddings shed some light on Problems \ref{prob:normal} and \ref{prob:dimm},
  whether the conclusion of Proposition \ref{prop:subspace} can be improved to show that if $(X,d)$
  is non-separable, then $\left(2^X, \mathscr T_{W(d)}\right)$ is non-normal is still unclear. This 
  leads to the following question. 

  \begin{question} \label{ques:closedcopy}
  Let $(X,d)$ be a non-separable metric space. Must $\left(2^X, \mathscr T_{W(d)}\right)$ contain a 
  closed copy of $\omega_1 \times (\omega_1+1)$?
  \end{question}

  In a recent paper \cite{Hernandez:2015}, Hern\'{a}ndez-Guti\'{e}rrez and Szeptycki also considered
  Problems \ref{prob:normal} and \ref{prob:dimm}. They proved that if $(X,d)$ is a locally separable
  metric space whose weight is a regular uncountable cardinal, then $\left(2^X, \mathscr T_{W(d)}
  \right)$ is non-normal. This result also answers partially to Question~\ref{ques:uniformdiscrete}. 
  Note that Question \ref{ques:closedcopy} also suggests the following relevant question:
  
  \begin{question}[\cite{cao-junnila:2014}] \label{ques:closedcopy2}
  Let $(X,d)$ be a metric space. If $\left(2^X, \mathscr T_{W(d)}\right)$ is non-normal, does it 
  contain a closed copy of $\omega_1 \times (\omega_1+1)$?
  \end{question}
  
  In \cite{hola:2015}, Hol\'{a} gave a partial answer to Question \ref{ques:closedcopy2}. In fact,
  she proved that the answer to Question \ref{ques:closedcopy2} is affirmative when every closed proper
  ball in $X$ is totally bounded. Consequently, the answer to Problems \ref{prob:normal} and \ref{prob:dimm}
  is also affirmative under this assumption.
  
  \section{A Problem on Completeness-type Properties\\ and Associated Questions}

  The study of completeness-type properties of Wijsman hyperspaces can be tracked back to Effros 
  \cite{effros:65}, whose main result can be interpreted as that a Polish space admits a metric whose 
  Wijsman topology is Polish. Beer \cite{beer:91} showed that the Wijsman hyperspace of any separable 
  complete metric space is Polish, and asked whether the Wijsman topology corresponding to an arbitrary 
  compatible metric for a Polish space is necessary Polish. Costantini \cite{costantini:95} answered
  affirmatively this problem, and a simpler proof of Costantini's theorem was given by Zsilinszky 
  \cite{zsilinszky:98} in terms of Choquet games. To the author's knowledge, the following problem was
  due to Beer in an oral communication.

  \begin{problem} \label{ques:beeroral}
  Is complete metrizability of $(X,d)$ (without separability) equivalent to any completeness-type 
  property of $\left(2^X, \mathscr T_{W(d)}\right)$?	
  \end{problem}

  Note that if $(X,d)$ is a non-separable metric space, the Wijsman topology $\mathscr T_{W(d)}$ 
  is Tychonoff but non-metrizable. Costantini \cite{costantini:98} showed that \v{C}ech-completeness 
  is not the right choice to answer Problem \ref{ques:beeroral}. Indeed, Costantini constructed a
  3-valued metric space on the set of real numbers whose Wijsman hyperspace is not \v{C}ech-complete. 
  More generally, in a recent paper \cite{cao-junnila:12}, Cao et al. established the following 
  embedding result.
  
  \begin{theorem}[\cite{cao-junnila:12}] \label{thm:universal}
  Every Tychonoff space can be embedded as a closed subspace in the Wijsman hyperspace of a complete 
  metric space $(X,d)$ which is locally $\mathbb R$.
  \end{theorem}
  
  In the light of Theorem \ref{thm:universal}, in addition to \v{C}ech-completeness, any 
  completeness-type property which is closed hereditary, is not the right choice to answer Problem 
  \ref{ques:beeroral}. This means that, to answer this problem, one has to turn attentions to 
  those completeness-type properties that are not closed-hereditary. Recall that a topological space 
  $X$ is said to be (resp. \emph{countably}) \emph{base compact} with respect to an open base 
  $\mathfrak B$ if $X$ is regular such that $\bigcap_{F \in \mathscr F} \overline{F} \ne \varnothing$ 
  for each (resp. countable) centered family ${\mathscr F} \subseteq \mathfrak B$, and $X$ is said to 
  be (resp. \emph{countably}) \emph{subcompact} with respect to an open base $\mathfrak B$ if $X$ is 
  regular such that $\bigcap_{F\in\mathscr F} \overline{F}\ne\varnothing$ for each (resp. countable) 
  regular filterbase ${\mathscr F}\subseteq\mathfrak B$. If ``regular" is replaced by ``quasi-regular" 
  and ``base" is replaced by ``$\pi$-base", the resulting spaces are called \emph{almost (countably)
  base compact, almost (countably) subcompact}, respectively. In literature, these properties are 
  called \emph{Amsterdam} properties. For details, refer to \cite{aarts:74}.
  
  \medskip
  Cao and Junnila \cite{cao-junnila:10} tackled Problem \ref{ques:beeroral} by considering the 
  Amsterdam properties of Wijsman hyperspaces. They discovered some interesting but peculiar results, 
  controversial to the results for other types of hyperspaces in \cite{cao:10}, \cite{cao-gutev:07} 
  and \cite{cao-tomita:07}.  
 
  \begin{example}[\cite{cao-junnila:10}] \label{exam:amsterdam1}
  There exists a metric space $(X, d)$ of the first category such that $\left(2^X, \mathscr T_{W(d)}
  \right)$ is countably base compact with respect to an open base $\mathfrak B$.	
  \end{example}
  
  \begin{example}[\cite{cao-junnila:10}] \label{exam:amsterdam2}
  There exists a separable metric space $(X, d)$ of the first category such that $\left(2^X, \mathscr 
  T_{W(d)}\right)$ is almost countably subcompact with respect to an open $\pi$-base $\mathfrak P$.
  \end{example}
  
  Next, we turn our attentions to some completeness-type properties defined by topological games. Let 
  $X$ be a topological space, and $\mathfrak P$ be a fixed open $\pi$-base. The \emph{Banach-Mazur game} 
  $BM(X)$ is played as follows: Players $\beta$ and $\alpha$ alternate in choosing elements of 
  $\mathfrak P$, with $\beta$ choosing first, so that 
  \[
  B_0 \supseteq A_0 \supseteq B_1 \supseteq A_1 \supseteq \cdots 
  \supseteq B_n \supseteq A_n \supseteq \cdots.
  \] 
  Then $B_0, A_0,\dots, B_n, A_n, \dots$ is a \emph{play} in $BM(X)$, and $\alpha$ wins this play if 
  $\bigcap_{n \in \mathbb N} A_n (= \bigcap_{n \in \mathbb N} B_n) \ne  \varnothing$, otherwise, $\beta$ 
  wins. A \emph{strategy} in $BM(X)$ is a function $\sigma: {\mathfrak P}^{< \omega} \to \mathfrak P$ 
  such that $\sigma(W_0, \dots, W_n) \subseteq W_n$ for all $n \in \mathbb N$, and $(W_0, \dots, W_n) 
  \in {\mathfrak P}^{n+1}$. A \emph{tactic} in $BM(X)$ is a function $t: {\mathfrak P} \to {\mathfrak 
  P}$ such that $t(W) \subseteq W$ for all $W \in {\mathfrak P}$. A \emph{winning strategy} (resp. 
  \emph{tactic}) for $\alpha$ is a strategy (resp. tactic) $\sigma$ such that $\alpha$ wins every play of 
  $BM(X)$ compatible with $\sigma$, i.e., such that $\sigma(B_0,\dots, B_n) = A_n$ (resp. $\sigma (B_n) 
  = A_n$) for all $n \in \mathbb N$. A winning strategy (resp. tactic) for $\beta$ is defined analogously. 
  The space $X$ is called (resp. \emph{weakly}) \emph{$\alpha$-favorable} \cite{telgarsky:87}, if $\alpha$ 
  has a winning tactic (resp. strategy) in $BM(X)$. The space $X$ is called \emph{$\beta$-favorable}, 
  if $\beta$ has a winning strategy in $BM(X)$. Let $\mathfrak B$ be an open base for $X$, and denote
  \[
  \mathscr E = \{ (x,U) \in X \times \mathfrak B: x \in U \}.
  \]
  The \emph{strong Choquet game} $Ch(X)$ is played similarly to the Banach-Mazur game. More precisely, 
  players $\beta$ and $\alpha$ alternate in choosing $(x_n, B_n) \in \mathscr E$ and $A_n \in \mathfrak B$, respectively, with $\beta$ choosing first so that for each $n \in \mathbb N$, $x_n \in A_n \subseteq B_n$, 
  and $B_{n+1} \subseteq A_n$. The play $(x_0, B_0), A_0, \dots, (x_n, B_n), A_n, \dots$ is won by $\alpha$, 
  if $\bigcap_{n\in \mathbb N} A_n (=\bigcap_{n \in \mathbb N} B_n) \ne \varnothing$; otherwise, $\beta$ 
  wins. A \emph{strategy} in $Ch(X)$ for $\alpha$ is a function $\sigma: {\mathscr E}^{<\omega} \to 
  \mathfrak B$ such that $x_n \in \sigma ((x_0, B_0),\dots, (x_n, B_n)) \subseteq B_n$ for all $((x_0, 
  B_0),\dots, (x_n, B_n)) \in {\mathscr E}^{<\omega}$. A \emph{tactic} in $Ch(X)$ for $\alpha$ is a 
  function $t : \mathscr E \to \mathfrak B$ such that $x \in t(x, B) \subseteq B$ for all $(x, B) \in 
  \mathscr E$. Winning strategies and tactics in $Ch(X)$ are defined similarly to the ones for the 
  Banach-Mazur game. The space $X$ is \emph{strongly $\alpha$-favorable} \cite{telgarsky:87} (resp. 
  \emph{strongly Choquet} \cite{kechris:94}), provided that $\alpha$ has a winning tactic (resp. 
  strategy) in $Ch(X)$.
  
  \medskip
  Regarding the completeness-type properties defined in the above, Pi\c{a}tkiewicz and Zsilinszky 
  \cite{piatkiewicz:10} established the following results. 

  \begin{theorem}[\cite{piatkiewicz:10}] \label{thm:stronglyalpha}
  Let $X$ be a locally separable metrizable space. The following are equivalent.
  \begin{itemize}
  \item[(1)] $\left(2^X, \mathscr T_{W(d)}\right)$ is strongly $\alpha$-favorable for every compatible 
  metric $d$ on $X$.
  \item[(2)] $\left(2^X, \mathscr T_{W(d)}\right)$  is strongly Choquet for every compatible metric $d$ 
  on $X$.
  \item[(3)] $X$ is completely metrizable.
  \end{itemize}
  \end{theorem}

  \begin{theorem}[\cite{piatkiewicz:10}] \label{thm:alpha}
  If $X$ is (weakly) $\alpha$-favorable and metrizable, then $\left(2^X, \mathscr T_{W(d)} \right)$ is 
  $\alpha$-favorable for every compatible metric $d$ on $X$.
  \end{theorem}

  \begin{example}[\cite{piatkiewicz:10}] \label{exam:firstcategory}
  There is a separable metric space $(X,d)$ of the first category such that $\left(2^X, \mathscr T_{W(d)}
  \right)$ is $\alpha$-favorable.
  \end{example}
   
  \begin{example}[\cite{piatkiewicz:10}] \label{exam:firstcategory2}
  There is a non-separable metric space $(X,d)$ of the first category such that $\left(2^X, \mathscr 
  T_{W(d)}\right)$  is strongly $\alpha$-favorable.
  \end{example}
   
  Note that Theorem \ref{thm:stronglyalpha} provides an answer to Problem \ref{ques:beeroral} in the 
  realm of locally separable metrizable spaces. In the light of Theorems \ref{thm:stronglyalpha} and 
  \ref{thm:alpha} as well as Examples \ref{exam:firstcategory} and \ref{exam:firstcategory2}, the 
  following two open questions are interesting.

  \begin{question}
  Let $X$ be a completely metrizable space. Must $\left(2^X, \mathscr T_{W(d)}\right)$ be strongly 
  Choquet or strongly $\alpha$-favorable for every compatible metric $d$ on $X$?
  \end{question}

  \begin{question}
  Let $X$ be a metrizable space. If $\left(2^X, \mathscr T_{W(d)}\right)$ is strongly Choquet or strongly 
  $\alpha$-favorable for every compatible metric $d$ on $X$, must $X$ be completely metrizable?
  \end{question}
  
  Recall that a topological space $X$ is said to be \emph{pseudocomplete} \cite{aarts:74} if $X$ is 
  quasi-regular and has a sequence $\{\mathfrak P_n : n \in \mathbb N\}$ of open $\pi$-bases such 
  that $\bigcap_{n\in \mathbb N} V_n \ne \varnothing$, whenever $\overline{V_{n+1}} \subseteq V_n 
  \in \mathfrak P_n$ for each $n\in \mathbb N$. Clearly, every almost countably subcompact space is 
  pseudocomplete.

  \begin{question}[\cite{cao-junnila:10}]
  If $(X,d)$ is pseudocomplete (resp. subcompact, base-compact), must $\left(2^X, \mathscr T_{W(d)}
  \right)$ be pseudocomplete (resp. subcompact, base-compact)?
  \end{question}
  
  \section{The Baire Property of Wijsman Hyperspaces}
  
  Recall that a topological space $X$ is said to be \emph{Baire} if the intersection of every sequence 
  of dense open subsets of $X$ is still dense. Note that a closed subspace of a Baire space may not be 
  Baire. A space is called \emph{hereditarily Baire} if every non-empty closed subspace is Baire. For 
  an alternative definition of a Baire space, refer to \cite{{haworth:77}}. The third basic problem on
  the Wijsman topology concerns the Baire property of Wijsman hyperspaces.
  
  \begin{problem} \label{prob:baire}
  Let $(X,d)$ be a metric space. Find characterizations for $\left(2^X, \mathscr T_{W(d)}\right)$ to 
  be a Baire space, in terms of some properties of $(X,d)$.
  \end{problem}

  Although it is not easy to find some completeness-type property for the Wijsman hyperspace of a 
  completely metrizable metric space, the following positive result on the Baire property was 
  discovered by Zsilinszky in 1996.

  \begin{theorem}[\cite{zsilinszky:96}] \label{thm:zsilinszky96}
  Let $(X,d)$ be a complete metric space. Then $\left(2^X, \mathscr T_{W(d)}\right)$ is a Baire space.
  \end{theorem}

  Note that, by Theorem \ref{thm:chaber-pol02} and Theorem \ref{thm:universal}, one should not expect 
  that the conclusion of Theorem \ref{thm:zsilinszky96} can be strengthened to ``strongly Baire". 
  However, this result can be strengthened as follows.

  \begin{theorem}[\cite{cao:10}] \label{thm:product}
  Let $X$ be a metrizable space. If $X^{\aleph_0}$ is Baire, then $\left(2^X, \mathscr T_{W(d)}\right)$ 
  is Baire for every compatible metric $d$ on $X$.
  \end{theorem}
  
  In \cite{cao-tomita:07}, Cao and Tomita constructed an example of a metric space $(X, d)$ such 
  that $X^n$ is a Baire space for each $n \in \mathbb N$, but $\left(2^X, \mathscr T_V\right)$ is 
  not a Baire space. Furthermore, Examples \ref{exam:amsterdam1}, \ref{exam:amsterdam2}, 
  \ref{exam:firstcategory} or \ref{exam:firstcategory2} imply that there exists a metric space $(X,d)$ 
  such that $\left(2^X, \mathscr T_{W(d)}\right)$ is Baire, but $\left(2^X, \mathscr T_V\right)$ is not 
  Baire. These examples suggest that the following question is interesting.

  \begin{question}
  Let $(X, d)$ be a metric space. If $X^n$ is Baire for each $n\in \mathbb N$,  must $\left(2^X, 
  \mathscr T_{W(d)}\right)$ be Baire?
  \end{question}

  At the 10th Prague Topological Symposium in 2006, Zsilinszky posed the following five open 
  questions in his talk. 

  \begin{question} \label{ques:Baire}
  If $(X,d)$ is a Baire metric space, must $\left(2^X, \mathscr T_{W(d)}\right)$ be Baire?
  \end{question}

  \begin{question} \label{ques:hereBaire}
  If $(X,d)$ is a hereditarily Baire metric space, must $\left(2^X, \mathscr T_{W(d)}\right)$ be Baire?
  \end{question}

  \begin{question} \label{ques:compatibleimpliesBaire}
  Let $X$ be a metrizable space. If $\left(2^X, \mathscr T_{W(d)}\right)$ is a Baire space for every 
  compatible metric $d$ on $X$, must $X$ be Baire?
  \end{question}

  \begin{question} \label{ques:wijsmanimpliesvietoris}
  Let $X$ be a metrizable space. If $\left(2^X, \mathscr T_{W(d)}\right)$ is a Baire space for every 
  compatible metric $d$ on $X$, must $\left(2^X, \tau_V\right)$ be Baire?
  \end{question}

  \begin{question} \label{ques:vietorisimplieswijsman}
  Let $X$ be a metrizable space. If $\left(2^X, \mathscr T_V \right)$ is a Baire space, must $\left(2^X, 
  \mathscr T_{W(d)}\right)$ be Baire for every compatible metric $d$ on $X$?
  \end{question}
  
  Question \ref{ques:hereBaire} was completely solved by Cao and Tomita in \cite{cao-tomita:10}, and its 
  answer is affirmative. There are several partial affirmative answers to Question \ref{ques:Baire}. 
  Note that if a metrizable space $X$ belongs to any of the following class of spaces:\\
  - Baire spaces having a countable open $\pi$-base;\\
  - separable Baire spaces;\\
  - hereditarily Baire spaces;\\
  - Baire spaces having a countable-in-itself open $\pi$-base;\\
  - almost locally $uK$-$U$ Baire spaces;\\
  - almost locally separable Baire spaces;\\
  - weakly $\alpha$-favorable spaces;\\
  - \v{C}ech-complete spaces;\\
  - spaces with any of the (countable) Amsterdam properties;\\
  - pseudocomplete spaces,\\
  then $X^{\aleph_0}$ is Baire, and thus by Theorem \ref{thm:product}, the answer to Question 
  \ref{ques:Baire} is affirmative. For details, refer to \cite{cao:10}. As mentioned at the end of 
  \cite{cao-junnila:10}, it is not possible to use the ``barely Baire spaces'' of Fleissner and Kunen 
  in \cite{fleissner-kunen:78} as counterexamples to Question \ref{ques:Baire}.
  
  \medskip
  The answer to Question \ref{ques:compatibleimpliesBaire} is affirmative in the class of almost locally 
  separable metrizable space, as the following result of Zsilinszky in \cite{zsilinszky:07} shows.
  
  \begin{theorem}[\cite{zsilinszky:07}] \label{thm:zsilinszky07}
  Let $X$ be an almost locally separable  metrizable space. Then $\left(2^X, \mathscr T_{W(d)}\right)$ is 
  Baire for every compatible metric $d$ on $X$ if and only if $X$ is Baire.
  \end{theorem}
  
  Note that Theorem \ref{thm:zsilinszky07} also provides a solution to Problem \ref{prob:baire} in the
  realm of almost locally separable metric spaces. However, the author does not know any information 
  toward (partial) solutions to Questions \ref{ques:wijsmanimpliesvietoris} and 
  \ref{ques:vietorisimplieswijsman}.


   \end{document}